\documentclass[10pt]{article}
\usepackage{amsfonts}
\usepackage{amsmath,hhline,epsfig,latexsym}
\usepackage{amssymb}
\usepackage[svgnames,dvipsnames,x11names]{xcolor} 
\usepackage{tikz}
\usepackage{graphicx}
\usepackage[english]{babel}
\usepackage{hyperref}
\usepackage{color}
\makeatletter
\renewcommand{\paragraph}{\@startsection{paragraph}{4}{0ex}%
   {-3.25ex plus -1ex minus -0.2ex}%
   {1.5ex plus 0.2ex}%
   {\normalfont\normalsize\bfseries}}
\makeatother

\stepcounter{secnumdepth}
\stepcounter{tocdepth}
\usepackage{amsmath}
\providecommand{\U}[1]{\protect\rule{.1in}{.1in}}
\newenvironment{claim}[1]{\par\noindent\underline{\bf Claim:}\space#1}{}
\newtheorem{thm}{Theorem}

\newtheorem{coro}{Corollary}

\newenvironment{proof}[1][Proof]
{\noindent\textbf{#1:} }{\hfill\rule{0.5em}{0.5em}}

\newcommand{\R}{\mathbb{R}}

\newcommand{\dis}{\displaystyle}

\newtheorem{theorem}{Theorem}[section]
\newtheorem{remark}[theorem]{Remark}

\title{\textbf{Non null controllability of Stokes equations with memory}}

\author{
	\textsc{Enrique Fern\' andez-Cara}\thanks{University of Sevilla, Dpto. E.D.A.N, Aptdo 1160, 41080 Sevilla, Spain. E-mail: {\tt cara@us.es}. Partially supported by grant MTM$2016$-$76990$-P, Ministry of Economy and Competitiveness (Spain).},
	\and
\textsc{Jos\'e Lucas F. Machado}\thanks{Department of Mathematics, Federal University of Pernambuco, UFPE, CEP 50740-545, Recife,
PE, Brazil. E-mail: {\tt  lucasmachado@dmat.ufpe.br}. Partially supported
by CNPq (Brazil).}  
	,
	\and
	\textsc{Diego A. Souza}\thanks{Department of Mathematics, Federal University of Pernambuco, UFPE, CEP 50740-545, Recife,
PE, Brazil. E-mail: {\tt diego.souza@dmat.ufpe.br}. Partially supported
by CNPq (Brazil) by the grants $313148$/$2017$-1 and Propesq (UFPE)-
Edital Qualis A. }
}
\date{}

\begin{document}
\maketitle



\begin{abstract}   
	In this paper, we study the null controllability of the three-dimensional Stokes equations with a memory term. For any positive final time $T>0$, we construct initial conditions
	such that the null controllability does not hold even if the controls act on the whole boundary. Moreover, we also prove that this negative result holds for distributed controls.
\end{abstract} 

%

\



\section{Introduction}

	Let $\Omega\subset\R^3$ be a smooth bounded domain and let $T > 0$ be a prescribed final time. 
	Let us introduce the Hilbert spaces
\[
	H(\Omega):=\{\,w\in L^2(\Omega)^3:\nabla\cdot w=0~\hbox{in}~\Omega,~w\cdot n=0~\hbox{on}~\partial\Omega\,\}
\]
	and
\[
	V(\Omega):=\{\,w\in H_0^1(\Omega)^3:\nabla\cdot w=0~\hbox{in}~\Omega\,\},
\]
	where $n=n(x)$ is the outward unit normal vector at $x\in\partial\Omega$. It is well know that $V(\Omega)\hookrightarrow H(\Omega)$ with a compact and dense embedding. 
	Consequently, identifying $H(\Omega)$ with its dual $H(\Omega)^\prime$, we have
\[
	V(\Omega)\hookrightarrow  H(\Omega)
	\hookrightarrow V(\Omega)^\prime,
\]
	with the second embedding being dense and compact.
		
		In the sequel, let us use the following notation\,: $Q:=\Omega\times (0,T)$ and $\Sigma:=\partial\Omega\times(0,T)$. The usual scalar product and norm in $L^2(\Omega)$ or $L^2(\Omega)^3$ 
		will be denoted by $(\cdot \,,\,\cdot)$ and $\|\,\cdot\,\|$, respectively. The symbols $C$, $C_0$, $C_1$, $\cdots$, will be used to design generic positive constants. 
		
		In this paper, let us consider the controlled Stokes equations with memory:
\begin{equation}\label{eq_stokes+memory}
\left\{
    \begin{array}{lcl}\dis
    		y_t-\Delta y-b\int_0^te^{-a(t-s)}\Delta y(\cdot\,,s) \ ds +\nabla p = 0 & \hbox{in} & Q,\\
    		  \noalign{\smallskip}\dis
    		\nabla\cdot y=0& \hbox{in} &Q,\\
    		y=v1_\gamma& \hbox{on} & \Sigma,\\
    		\noalign{\smallskip}\dis
    		y(\cdot\,,0)=y^0 &\hbox{in} &\Omega,
    \end{array}
    \right.
\end{equation}
	 where $a,\,b>0$ and $\gamma\subset \partial\Omega$ is a non-empty open subset of the boundary. 
	 Here, $v \in L^2(\gamma\times(0,T))$ is a control function which is acting on $\gamma$ during the whole interval $(0,T)$ and $y^0\in H(\Omega)$ is an initial data.

	{It is well known that for any $y^0\in H(\Omega)$ and $v\in L^2(\gamma\times(0,T))$, there exists exactly one solution to \eqref{eq_stokes+memory}, in the sense of transposition. 
	In other words, there exists a unique  $y\in L^2(0,T;H(\Omega))\cap C^0([0,T];V(\Omega)^\prime)$ satisfying  
	\begin{equation*}
		\int^T_0(y(\cdot\,,t),g(\cdot\,,t))\ dt=(y^0,\psi(\cdot\,,0))-\iint_{\gamma\times(0,T)}v\left(-\pi n+ \dfrac{\partial \psi}{\partial n} + b\int^T_t e^{-a(s-t)}\dfrac{\partial \psi}{\partial n}(\cdot\,,s) \, ds \right) \! d\Gamma dt,
	\end{equation*}
for all $g\in L^2(0,T;H(\Omega))$, where $\psi$ is, together with some pressure $\pi$, the unique (strong) solution to
\begin{equation*}
\left\{
    \begin{array}{lll}\dis
    		-\psi_t-\Delta \psi-b\int_t^Te^{-a(s-t)}\Delta \psi(\cdot\,,s) \ ds +\nabla \pi = g & \hbox{in} & Q,\\
    		  \noalign{\smallskip}\dis
    		\nabla\cdot \psi=0& \hbox{in} &Q,\\
    		\psi=0& \hbox{on} & \Sigma,\\
    		\noalign{\smallskip}\dis
    		\psi(\cdot\,, T)=0 &\hbox{in} &\Omega.
    \end{array}
    \right.
\end{equation*} 
	Also, it can be proved that $y\in C^0([0,T]; V_\sigma(\Omega)^\prime)$ for all $\sigma>1/2$, where $V_\sigma(\Omega)=H^\sigma(\Omega)\cap H(\Omega)$; for more details, see \cite{Lions_Magenes,Teman}. Of course, if $v1_\gamma$ is regular enough (for instance, $v = \overline{y}|_{\gamma\times(0,T)}$ with $\overline{y}\in L^2(0,T;V(\Omega))$ and $\overline{y}_t\in L^2(0,T;V(\Omega)^\prime)$), then $y$ is, together 
	with some pressure $p$, the unique weak solution to \eqref{eq_stokes+memory}.

	On the other hand, notice that $y$ is the unique function in $C^0([0,T];V(\Omega)^\prime)$ such that, for~all~$\overline{t}\in [0,T]$, one has
	\begin{equation*}
		\langle y(\cdot\,,\overline{t}),w\rangle =(y^0,\psi(\cdot\,,0))-\!\!\!\iint_{\gamma\times(0,\overline{t})}\!\!\!\!v\left(-\pi n+ \dfrac{\partial \psi}{\partial n} + b\int^T_{t} e^{-a(s-t)}
		\dfrac{\partial \psi}{\partial n}(\cdot\,,s) \, ds \right) \! d\Gamma dt
	\end{equation*}
	for all $w\in V(\Omega)$, where $(\psi, \pi)$ is the unique (strong) solution to 
\begin{equation*}
\left\{
    \begin{array}{lll}\dis
    		-\psi_t-\Delta \psi-b\int_t^Te^{-a(s-t)}\Delta \psi(\cdot\,,s) \ ds +\nabla  \pi = 0 & \hbox{in} & \Omega\times (0,\overline{t}),\\
    		  \noalign{\smallskip}\dis
    		\nabla\cdot \psi=0& \hbox{in} &\Omega\times (0,\overline{t}),\\
    		\psi=0& \hbox{on} & \partial\Omega\times (0,\overline{t}),\\
    		\noalign{\smallskip}\dis
    		\psi(\cdot\,, \overline{t})=w &\hbox{in} &\Omega.
    \end{array}
    \right.
\end{equation*} 
	}

	The boundary null controllability property for \eqref{eq_stokes+memory} is as follows: for each $y^0\in H(\Omega)$, find a boundary control $v\in L^2(\gamma\times(0,T))$ 
	such that the associated solution satisfies $y(\cdot,T)=0$.  
	
	When $b = 0$, \eqref{eq_stokes+memory} is the Stokes equations and it is well known that the null controllability holds. 
	In the general case, the presence of the memory term brings new difficulties in the analysis of the controllability for \eqref{eq_stokes+memory} .  

	By a duality argument, it is not difficult to see that the null controllability of \eqref{eq_stokes+memory} is equivalent to prove an observability inequality for the adjoint system:
\begin{equation}\label{1a}
\left\{
    \begin{array}{lll}\dis
    		-\varphi_t-\Delta \varphi-b\int_t^Te^{-a(s-t)}\Delta \varphi(\cdot\,,s) \ ds +\nabla q = 0 & \hbox{in} & Q,\\
    		  \noalign{\smallskip}\dis
    		\nabla\cdot \varphi=0& \hbox{in} &Q,\\
    		\varphi=0& \hbox{on} & \Sigma,\\
    		\noalign{\smallskip}\dis
    		\varphi(\cdot\,, T)=\varphi^0 &\hbox{in} &\Omega.
    \end{array}
    \right.
\end{equation} 
The usual way to deduce such an observability estimate is to first prove a global Carleman inequality. But it seems difficult to adapt this approach in the presence of an
integro-differential term.

	In the last decades, a lot of researchers has focused attention on the controllability of systems governed by linear and nonlinear PDEs. 
	For linear PDEs, the first main contributions were obtained in \cite{Imanuvilov,Lebeau,Lions1,Russell2,Russell}. 
	For instance, in  \cite{Russell}, D.L. Russell present a rather complete survey on the most relevant results available at that time
	and, in \cite{Lions1}, J.-L. Lions introduce the so called {\it Hilbert Uniqueness Method} (HUM for short). 
 For semilinear systems of PDEs, one can find the first contributions in \cite{Fabre, Fursikov,Triggiani,Zuazua} and  some other related results can be found in \cite{Coron,glowinski}.


	In the context of fluid mechanics, the main controllability results are related to the  Burgers, Stokes, Euler and Navier-Stokes equations. 
	For Stokes equations, the approximate and null controllability with distributed controls have been established in \cite{Fabre_stoke,Imanuvilov_stokes}, respectively. 
	For the Euler equations, global controllability results are proved in \cite{Coron_euler,Glass}. On the other hand, for the Navier-Stokes equations with initial and Dirichlet boundary conditions 
	only local controllability results are available, see, for instance, \cite{Cara_NS, Cara_n-1,Fursikov_NS, Imanuvilov_stokes}.


	For an one-dimensional heat equation with memory, the lack of null controllability for a large class of memory kernels and controls was established in \cite{Pandolfi}, where the notion of null controllability also 
	requeres that $\int_0^Ty(\cdot,t) \ dt=0$. In higher dimensional case, Guerrero and Imanuvilov proved, in \cite{Guerrero}, that 
	the null controllability does not hold for the following heat equation with memory:
\begin{equation}\label{eq_guerrero}
\left\{
    \begin{array}{lll}\dis
    		y_t-\Delta y-\int_0^t\Delta y(\cdot\,,s) \ ds= 0 & \hbox{in} & Q,\\
    		y=v1_\gamma & \hbox{on} & \Sigma,\\
    		\noalign{\smallskip}\dis
    		y(\cdot\,,0)=y^0 &\hbox{in} &\Omega.
    \end{array}
    \right.
\end{equation}
	{A similar result was obtained in \cite{Zhou_Gao} by Zhou and Gao for the system
\begin{equation*}\label{eq_Gao}
\left\{
    \begin{array}{lll}\dis
    		y_t-\Delta y-b\int_0^te^{-a(t-s)} y(\cdot\,,s) \ ds = 0 & \hbox{in} & Q,\\
    		y=v& \hbox{on} & \Sigma,\\
    		\noalign{\smallskip}\dis
    		y(\cdot\,,0)=y^0 &\hbox{in} &\Omega.
    \end{array}
    \right.
\end{equation*}}

	Our main goal in this work is to prove that the null controllability of  $(\ref{eq_stokes+memory})$ does not hold.  More precisely, we have the following result:
\begin{thm}\label{theo_main}
	Let $T>0$ be given. There exists initial data $y^0 \in H(\Omega)$ such that, for any control $v \in L^2(\gamma\times(0,T))$, the associated solution to \eqref{eq_stokes+memory} is not identically 
	zero at time $T$.
\end{thm}
	The proof of this theorem follows some ideas of~\cite{Guerrero}. Indeed, we prove that observability inequality does not hold and then we construct explicitly initial data that cannot 
	be steered to zero.

	We also have a negative result related to distributed control systems:
\begin{equation}\label{eq_stokes+memory_distributed}
\left\{
    \begin{array}{lll}\dis
    		y_t-\Delta y-b\int_0^te^{-a(t-s)}\Delta y(\cdot\,,s) \ ds +\nabla p = v1_\omega & \hbox{in} & Q,\\
    		  \noalign{\smallskip}\dis
    		\nabla\cdot y=0& \hbox{in} &Q,\\
    		y=0& \hbox{on} & \Sigma,\\
    		\noalign{\smallskip}\dis
    		y(\cdot\,,0)=y^0 &\hbox{in} &\Omega,
    \end{array}
    \right.
\end{equation}
	where $\omega\subset \Omega$ is an open subset. More precisely, as an immediate consequence of Theorem \ref{theo_main}, we have the following result:
\begin{coro}\label{coro_main}
	Let $T>0$ be given and $\omega\subset\Omega$ is a non-empty open proper subset. There exists initial data $y^0 \in H(\Omega)$ such that, for any control 
	$v \in L^2(\omega\times(0,T))$, the associated solution to \eqref{eq_stokes+memory_distributed} is not identically zero at time $T$.
\end{coro}

\begin{remark}
	 Theorem \ref{theo_main} and Corollary \ref{coro_main} also hold for the following Stokes equations with memory:
\begin{equation*}\label{eq_stokes+memory}
\left\{
    \begin{array}{lcl}\dis
    		y_t-\Delta y+b\int_0^te^{-a(t-s)} y(\cdot\,,s) \ ds +\nabla p = 0 & \hbox{in} & Q,\\
    		  \noalign{\smallskip}\dis
    		\nabla\cdot y=0& \hbox{in} &Q,\\
    		y=v1_\gamma& \hbox{on} & \Sigma,\\
    		\noalign{\smallskip}\dis
    		y(\cdot\,,0)=y^0 &\hbox{in} &\Omega.
    \end{array}
    \right.
\end{equation*}
\end{remark}

	The control analysis of \eqref{eq_stokes+memory} and \eqref{eq_stokes+memory_distributed} is motivated by the interest to understand the limits of controlling viscoelastic fluids of the Oldroyd's kind. 
	Specifically, let us consider the following boundary and distributed nonlinear  control systems:
\begin{equation}\label{old_2}
\left\{
    \begin{array}{lll}
    		y_t-\nu\Delta y+(y\cdot\nabla)y +\nabla p= \nabla \cdot\tau & \hbox{in} & Q,\\
    		  \noalign{\smallskip}\dis
    		 \tau_t+(y\cdot\nabla)\tau+g(\nabla y,\tau)+a\tau=2b D( y)& \hbox{in} & Q,\\
    		\noalign{\smallskip}\dis
    		\nabla\cdot y=0& \hbox{in} & Q,\\
   		\noalign{\smallskip}\dis
   		y=v 1_\gamma& \hbox{on} & \Sigma,\\
   		\noalign{\smallskip}\dis
  		y(\cdot\,,0)=y^0, \quad \tau(\cdot\,,0)=\tau^0 &\hbox{in} &\Omega
    \end{array}
    \right.
\end{equation}
and
\begin{equation}\label{old_1}
\left\{
    \begin{array}{lll}
    		y_t-\nu\Delta y+(y\cdot\nabla)y +\nabla p= \nabla \cdot\tau +v 1_\omega& \hbox{in} & Q,\\
    		\noalign{\smallskip}\dis
    		\tau_t+(y\cdot\nabla)\tau+g(\nabla y,\tau)+a\tau=2b D (y)& \hbox{in} & Q,\\
    		\noalign{\smallskip}\dis
    		\nabla\cdot y=0& \hbox{in} & Q,\\
   		\noalign{\smallskip}\dis
   		y=0& \hbox{on} & \Sigma,\\
   		\noalign{\smallskip}\dis
  		y(\cdot\,,0)=y^0, \quad \tau(\cdot\,,0)=\tau^0 &\hbox{in} &\Omega,
    \end{array}
    \right.
\end{equation}
	where 
	$g(\nabla y,\tau):=\tau W(y)-W(y ) \tau-k[D(y)\tau+\tau D(y)]$,  $k\in[-1,1]$, $D(y):={1\over 2}(\nabla y+\nabla y^t)$ and $W(y):=\frac{1}{2}(\nabla y-\nabla y^t)$.
	The functions $y$, $p$ and $\tau$ are the velocity field, the pressure distribution and the elastic extra-stress tensor of the fluid, respectively;
	$y^0 \in H(\Omega)$ and $\tau^0\in L^2(\Omega; \mathcal{L}_s(\R^3))\footnote{$\mathcal{L}_s(\R^3)$ is the space of symmetric real 3 $\times$ 3 matrices.}$.

	The theoretical analysis of the Oldroyd systems $(\ref{old_2})$ and $(\ref{old_1})$ has been the subject of considerable work. Notice that these systems  are more difficult to solve than the 
	usual Navier-Stokes equations. The main reason is the presence of the nonlinear term $g(\nabla y,\tau)$. For more details, see \cite{Enrique_oldroyd, Lions, Renardy_3}.

	It is worth to mentioning that in \cite{Doubova} the authors  studied a linear version of \eqref{old_1}:
\begin{equation}\label{eq_jeffreys}
\left\{
    \begin{array}{lll}
    		y_t-\Delta y +\nabla p= \nabla \cdot\tau +v 1_\omega& \hbox{in} & Q,\\
    		\noalign{\smallskip}\dis
    		\tau_t+a\tau=2b D(y)& \hbox{in} & Q,\\
    		\noalign{\smallskip}\dis
    		\nabla\cdot y=0& \hbox{in} & Q,\\
   		\noalign{\smallskip}\dis
   		y=0& \hbox{on} & \Sigma,\\
   		\noalign{\smallskip}\dis
  		y(\cdot\,,0)=y^0, \quad \tau(\cdot\,,0)=\tau^0 &\hbox{in} &\Omega.
    \end{array}
    \right.
\end{equation}
	Plugging the explicit solution $\tau$ of \eqref{eq_jeffreys}$_2$ in  \eqref{eq_jeffreys}$_1$, the system above  can be equivalently rewritten as an integro-differential equation 
	in $y$:
\begin{equation}\label{eq_jeff_int}
\left\{
    \begin{array}{lll}
    		\dis	 y_t-\Delta y -b\int_0^te^{-a(t-s)}\Delta y(\cdot\,,s) \ ds+\nabla p= e^{-at}\nabla \cdot\tau^0 +v 1_\omega& \hbox{in} & Q,\\
    		\noalign{\smallskip}\dis
    		\nabla\cdot y=0& \hbox{on} & Q,\\
    		\noalign{\smallskip}\dis
    		y=0& \hbox{on} & \Sigma,\\
   		\noalign{\smallskip}\dis
  		y(\cdot\,,0)=y^0 &\hbox{on} &\Omega.
    \end{array}
\right.
\end{equation}

	In \cite[Theorem $1.1$ and $1.2$]{Doubova}, approximate controllability results are established for \eqref{eq_jeffreys} in both distributed and boundary cases.
	Notice that, if $\tau^0$ is the null matrix then the system \eqref{eq_jeff_int} is exactly the system \eqref{eq_stokes+memory_distributed}.

	The linear system \eqref{eq_jeffreys} governs the behavior of linear  viscoelastic fluids of the Jeffreys' kind. This system without the viscosity term corresponds to the so called {\it linear Maxwell fluids}. In \cite{Boldrini},
	the authors have proved large time approximate-finite dimensional and exact controllability results for some suitable control domains;
	 for others works, see \cite{Renardy_1, Renardy_2}. On the physical meanings of these systems, see for instance \cite{Joseph, Renardy}.

	It is worth mentioning that in \cite{Doubova} the null controllability of linear Jeffreys fluids is formulated as an open problem.
	Theorem \ref{theo_main} and Corollary \ref{coro_main} solve this open question proving that the null controllability does not hold.

	This paper is organized as follows. In Section \ref{section_eigenfunction}, we compute the eigenfunctions and eigenvalues of the  Stokes operator in a ball and we prove some relevant estimates. 
	In Section \ref{section_demonstration}, we prove Theorem \ref{theo_main}. Finally, in Section \ref{coments}, we present 
	 some additional comments and open problems.

%

\section{The radially symmetric eigenfunctions of the Stokes operator}\label{section_eigenfunction}

	In this section, let us assume that $\Omega$ is the ball of radius $R$ and centered at the origin. Let us compute explicitly the eigenfunctions and eigenvalues of the Stokes operator and then we will 
	deduce some crucial estimates that will be used to prove Theorem \ref{theo_main}. For simplicity, the coordinates of a generic point in $\Omega$ will be denoted by $x$, $y$ and $z$.

	Let us compute $(\varphi, q)$ and $\lambda$ such that
\begin{equation}\label{prob_autovalor}
\left\{
    \begin{array}{lll}
    		\dis -\Delta \varphi +\nabla q= \lambda\varphi& \hbox{in} & \Omega,\\
    		\noalign{\smallskip}\dis
    		\nabla\cdot \varphi=0& \hbox{in} & \Omega,\\
    		\noalign{\smallskip}\dis
    		\varphi=0& \hbox{on} & \partial\Omega.
    \end{array}
\right.
\end{equation}

	Let us looking for eigenfunctions as the curl of radial stream functions, i.e. $\varphi=\nabla\times \psi$, for some radial stream function $\psi$. Setting $w=\nabla\times\varphi$, we can easily deduce that if
	$(w,\psi)$ solves following
	 the eigenvalue problem
\begin{equation}\label{eq:wpsi}
\left\{
    \begin{array}{lll}
    		\dis -rw^{\prime\prime}-2w^{\prime}= \lambda rw& \hbox{in} & (0,R),\\
    		\noalign{\smallskip}\dis
    		-r\psi^{\prime\prime}-2\psi^{\prime}= rw& \hbox{in} & (0,R),\\
    		\noalign{\smallskip}\dis
    		\psi(R)=0, \ \ \psi^{\prime}(R)=0, \ \ \lambda>0&  & 
    \end{array}
\right.
\end{equation}
	then $\varphi=\nabla \times \psi$ is a solution for the eigenvalue problem \eqref{prob_autovalor}. Here, we are using the notation $r=\sqrt{x^2+y^2+z^2}$ for any $(x,y,z)\in \Omega$.

	In order to compute the solution of the previous problem, let us make the following change of variables: $\zeta=rw$ and $\phi=r\psi$. Then, from \eqref{eq:wpsi}, we see that $\zeta$ and $\phi$ satisfy
\begin{equation*}
\left\{
    \begin{array}{lll}
    		\dis -\zeta^{\prime\prime}= \lambda \zeta, \quad -\phi^{\prime\prime}= \zeta& \hbox{in} & (0,R),\\
    		\noalign{\smallskip}\dis
    		\zeta(0)=0, \ \ \phi(0)=0,\\
    		\noalign{\smallskip}\dis
    		\phi(R)=0, \ \ \phi^{\prime}(R)=0, \ \ \lambda>0.&  & 
    \end{array}
\right.
\end{equation*}

This way, it is not difficult compute explicitly the eigenvalues $\lambda_n$ and the corresponding eigenfunctions $(\varphi_n,q_n)$ for \eqref{prob_autovalor}:
\begin{equation}\label{def_autofunc}
\left\{
\begin{array}{l}
	\noalign{\smallskip}\dis
	\varphi_n(x,y,z)=\dfrac{1}{\lambda_n^{1/2}r^2}\left(\cos(\lambda_n^{1/2}r)-\dfrac{1}{\lambda_n^{1/2}r}\sin(\lambda_n^{1/2}r)\right)(y-z,z-x,x-y),\\
	\noalign{\smallskip}\dis
	q_n\equiv0,\\
	\noalign{\smallskip}\dis
	\lambda^{1/2}_nR=tg(\lambda^{1/2}_nR).
\end{array}
\right.
\end{equation}

	Notice that
\begin{equation}\label{prop_lambda}
\lambda_n=\frac{\pi^2}{R^2}(n+1/2)^2-\varepsilon_n, \quad \hbox{for some}\quad\varepsilon_n>0 \quad \hbox{with}\quad \varepsilon_n \to  0.
\end{equation}

	It is not difficult to see that $\{\varphi_n\}_{n\in\mathbb{N}}$ is an orthogonal family in $H(\Omega)$. Also, using \eqref{def_autofunc}$_3$, we can compute the $L^2$-norm of $\varphi_n$:
\begin{equation}\label{est_autofuncao_1}
\begin{alignedat}{2}
    		\dis \|\varphi_n\|^2&=\dis 8\pi\int_0^R\left(\dfrac{\cos(\lambda_n^{1/2}r)}{\lambda_n^{1/2}}-\dfrac{\sin(\lambda_n^{1/2}r)}{\lambda_nr}\right)^2dr  \\
    		\noalign{\smallskip}\dis
    		&=\dfrac{8\pi}{\lambda_n^{3/2}}\left[ \dfrac{\lambda_n^{1/2}R}{2}+\sin(\lambda_n^{1/2}R)\left( \dfrac{\cos(\lambda_n^{1/2}R)}{2}-\dfrac{\sin(\lambda_n^{1/2}R)}{\lambda_n^{1/2}R}  \right)   \right] \\
    		\noalign{\smallskip}\dis
    		&= \dfrac{2\pi R}{\lambda_n}(1-\cos(2\lambda_n^{1/2}R)).
    \end{alignedat}
\end{equation}

	Now, from \eqref{prop_lambda} and \eqref{est_autofuncao_1}, if $n$ is large enough then we obtain that~$\cos(2\lambda_n^{1/2}R)
<0$ and, consequently
\begin{equation}\label{est_autofuncao_2}
\|\varphi_n\|^2\geq \dfrac{2\pi R}{\lambda_n}.
\end{equation}

	Finally, we will deduce some estimates for the normal derivatives of $\varphi_n$. Using \eqref{def_autofunc}$_1$ and \eqref{def_autofunc}$_3$, we easily compute:
\begin{equation*}
	\begin{alignedat}{2}
	\noalign{\smallskip}\dis
		\dfrac{\partial \varphi^1_n}{\partial n}\biggr|_{\partial\Omega}&=\dis\biggl( -\dfrac{\sin(\lambda_n^{1/2}R)}{R^2}-3\dfrac{\cos(\lambda_n^{1/2}R)}{\lambda^{1/2}_nR^3}+ 3\dfrac{\sin(\lambda_n^{1/2}R)}{\lambda_nR^4}\biggr)(y-z),\\
		\noalign{\smallskip}\dis
		\dfrac{\partial \varphi^2_n}{\partial n}\biggr|_{\partial\Omega}&=\dis\biggl( -\dfrac{\sin(\lambda_n^{1/2}R)}{R^2}-3\dfrac{\cos(\lambda_n^{1/2}R)}{\lambda^{1/2}_nR^3}+ 3\dfrac{\sin(\lambda_n^{1/2}R)}{\lambda_nR^4}\biggr)(z-x),\\
		\noalign{\smallskip}\dis
		\dfrac{\partial \varphi^3_n}{\partial n}\biggr|_{\partial\Omega}&=\dis\biggl( -\dfrac{\sin(\lambda_n^{1/2}R)}{R^2}-3\dfrac{\cos(\lambda_n^{1/2}R)}{\lambda^{1/2}_nR^3}+ 3\dfrac{\sin(\lambda_n^{1/2}R)}{\lambda_nR^4}\biggr)(x-y).
	\end{alignedat}
\end{equation*}
Thanks to  \eqref{def_autofunc}$_3$, we get the relation
\begin{equation*}
	\begin{array}{lll}
	 -3\dfrac{\cos(\lambda_n^{1/2}R)}{\lambda^{1/2}_nR^3}+ 3\dfrac{\sin(\lambda_n^{1/2}R)}{\lambda_nR^4}
	 &=&0,
 	\end{array}
\end{equation*}
	and then:
\begin{equation}\label{derivada_normal_varphi_n}
\!\!\!\!\dfrac{\partial \varphi_n}{\partial n}\biggr|_{\partial\Omega}\!=\!-\dfrac{\sin(\lambda_n^{1/2}R)}{R^2}(y-z,z-x,x-y).
\end{equation}

%
%

\section{Lack of null controllability}\label{section_demonstration}

 In this section, we prove Theorem \ref{eq_stokes+memory}. We will follow some ideas presented in \cite{Guerrero}.
 
 \begin{remark}
 	Notice that in order to prove of Theorem \ref{theo_main} it is sufficient to consider the case where $\Omega$ is a ball and the solution is supposed to be radially symmetric. 
	Indeed, if­ $\Omega$ is a general bounded domain in $\R^3$, we consider an open ball $B\subset \Omega$. Once the result is established for the domain $B$, we have that, for any positive $T$, there exists 
	an initial condition $\widehat{y}^0\in H(B)$ such that, for any boundary control $v\in L^2(\partial B\times (0,T))$, the associated solution $\hat{y}$ is not identically equal to zero at time $T$. 
	Now, let us extend (by zero)  to $\Omega$ the initial condition $\widehat{y}^0$ and consider the system \eqref{eq_stokes+memory} in $\Omega$. Therefore, arguing by contradiction,
	 one can easily verify that for this extended initial condition the null controllability for \eqref{eq_stokes+memory} at time $T$ also fails. 
 \end{remark}
 
	From the previous remark, let us consider the case that $\Omega$ is a ball of radius $R$. 
%
  It is well known that the null controllability
of \eqref{eq_stokes+memory} is equivalent to the following observability
inequality for the solutions to \eqref{1a}: 
 \begin{equation}\label{des-obs}
 \|\varphi(\,\cdot\,,0)\|^2 \leq C\!\iint_\Sigma \left|\left(-q\text{Id}+\nabla\varphi+b\!\int_t^T\!\!\!e^{-a(s-t)}\nabla \varphi(\cdot,s) \ ds\right)\cdot n\right|^2d\Gamma dt \quad  \forall \varphi^0\in H(\Omega).
 \end{equation}
 
 	The goal is to show that there is no positive constant $C$ such that \eqref{des-obs} is true. To this purpose, we will construct a family of solutions $\varphi^M$ to \eqref{1a}, for all sufficiently large $M$,
	 such that 
 \begin{equation}\label{eq:phiM}
 	\|\varphi^M(\,\cdot\,,0)\|\geq \dfrac{C_1}{M^6} 
 \end{equation}
and 
  \begin{equation}\label{eq:normalphiM}
 	\iint_\Sigma \left|\left(-q\text{Id}+\nabla\varphi^M+b\int_t^Te^{-a(s-t)}\nabla \varphi^M(\cdot,s) \ ds\right)\cdot n\right|^2d\Gamma dt\leq \dfrac{C_2}{M^{10}}, 
 \end{equation}
 for some 
 positive constants $C_1$ and $C_2$ (independent of $M$).
 
  	Therefore, using these properties of $\varphi^M$, we will be able to construct initial conditions ~$\overline{y}^0$ in $H(\Omega)$ such that 
	the solution to \eqref{eq_stokes+memory} cannot be steered to zero, no matter the control is.

\subsection{Construction of $\varphi^M$}\label{sec:3.1}

	Let us now present some computation which will inspire the construction of $\varphi^M$. To do this, let us first consider
\begin{equation*}
\varphi^0=\sum_{n\geq 1}\beta_n\varphi_n, \quad \{\beta_n\}\in \ell ^2.
\end{equation*}
Then, the associated solution to \eqref{1a} can be written in the form
 \begin{equation}\label{def_varphi}
\varphi(\,\cdot\,,t)=\sum_{n\geq 1}\alpha_n(t)\varphi_n\quad \forall  t \in (0,T),
\end{equation}
 where the functions $\alpha_n$ satisfy the following second-order Cauchy problem:
\begin{equation}\label{EDO_alpha}
\left\{
\begin{array}{lll}
  	-\alpha_n^{\prime\prime}+(\lambda_n+a)\alpha_n^{\prime}-\lambda_n(a+b)\alpha_n =0 & \hbox{in} & (0,T),\\
    		  \noalign{\smallskip}\dis
    		\alpha_n(T)=\beta_n,&&\\
    		  \noalign{\smallskip}\dis
	  \alpha_n^{\prime}(T)=\lambda_n\beta_n.&&
    \end{array}
    \right.
\end{equation} 
	\begin{remark}\label{rmk:q}
		In view of \eqref{def_autofunc}$_2$, we have that the pressure $q$ associated to $\varphi$ in \eqref{def_varphi} is zero.
	\end{remark}

 	It is clear that, there exists $n_0\in \mathbb{N}$~such that if $n\geq n_0$ then one has $D_n:=(\lambda_n+a)^2-4(a+b)\lambda_n>0$. This way, taking 
	$\beta_n=0$ for $n<n_0$, we have that
 \begin{equation}\label{def_alpha_n}
 \left\{
\begin{array}{lll}
\alpha_n\equiv0\quad \forall n< n_0,\\
    		  \noalign{\smallskip}\dis
 \alpha_n(t)=C_{1,n}e^{\mu_n^+(T-t)}+C_{2,n}e^{\mu_n^-(T-t)}\quad \forall t\in (0,T),\quad\forall n\geq n_0,
    \end{array}
    \right.
 \end{equation}
 where
 \begin{equation}\label{def_mu}
 \mu_n^+=-\dfrac{(\lambda_n+a)+\sqrt{D_n}}{2}\quad\hbox{and}\quad \mu_n^-=-\dfrac{(\lambda_n+a)-\sqrt{D_n}}{2}
 \end{equation} 
 and the constants $C_{1,n}$ and $C_{2,n}$ are given by
 \begin{equation}\label{def_C_n}
 C_{1,n}
 =\beta_n	\dfrac{\lambda_n-a+\sqrt{D_n}}{2\sqrt{D_n}}
 \quad \text{and} \quad
 C_{2,n}
 =\beta_n\dfrac{a-\lambda_n+\sqrt{D_n}}{2\sqrt{D_n}}.
 \end{equation}

\begin{remark}
	It is not difficult to see that  $\mu_n^+\to-\infty$ and $\mu_n^-\to -(a+b)$ as $n\to+\infty$.
\end{remark}
 
%

 	Using \eqref{est_autofuncao_2}, \eqref{def_varphi}, \eqref{def_alpha_n} and the orthogonality of $\varphi_n$, we see that
 \begin{equation}\label{brasil}
 	\begin{alignedat}{2}\dis
    		\|\varphi(\,\cdot\,,0)\|^2& =&
		\dis\sum_{n\geq n_0}(C_{1,n}e^{\mu_n^+T}+C_{2,n}e^{\mu_n^-T})^2\|\varphi_n\|^2\\
    		  \noalign{\smallskip}\dis
    		  &\geq&\dis\sum_{n\geq n_0}\dfrac{2\pi R}{\lambda_n}(C_{1,n}e^{\mu_n^+T}+C_{2,n}e^{\mu_n^-T})^2.
    \end{alignedat}
 \end{equation}

	In order to estimate the right hand side of \eqref{des-obs} from above, it is sufficient to find an estimate for the term
\[
	\iint_\Sigma \left|\dfrac{\partial \varphi}{\partial n}\right|^2d\Gamma dt.
\]
	In order to simplify the computations, let us introduce the weight $e^{2(a+b)(T-t)}$ in the above integral and estimate it, i.e. let us bound the term
\begin{equation*}\label{alemanha}
\iint_\Sigma e^{2(a+b)(T-t)}\left|\dfrac{\partial \varphi}{\partial n}\right|^2d\Gamma dt.
\end{equation*}

 Taking into account \eqref{derivada_normal_varphi_n}, the following estimate holds:
\begin{equation*}
	\begin{alignedat}{2}
		\left|\dfrac{\partial \varphi}{\partial n}\right|^2
		\leq\dis12\left|\sum_{n\geq n_0}\gamma_n\alpha_n(t)\right|^2,
	\end{alignedat}
\end{equation*}
 	where $\gamma_n:=\sin(\lambda_n^{1/2}R)/R$. Therefore, 
\begin{equation}\label{japao}
	\begin{alignedat}{2}
	\dis\iint_\Sigma e^{2(a+b)(T-t)}\left|\dfrac{\partial\varphi}{\partial n}\right|^2d\Gamma dt\leq&~ 
	 \dis48\pi R^2\int_0^Te^{2(a+b)(T-t)}\left|\sum_{n\geq n_0}\alpha_n(t)\gamma_n\right|^2dt\\
	 	\noalign{\smallskip}\dis
	 \leq&~\dis96\pi R^2\int_0^T\left(\sum_{n\geq n_0}\gamma_n C_{1,n}e^{(\mu_n^++a+b)(T-t)}   \right)^2dt\\
	 	\noalign{\smallskip}\dis
	 &\dis+\ 96\pi R^2\int_0^T\left(\sum_{n\geq n_0}\gamma_nC_{2,n}e^{(\mu_n^-+a+b)(T-t)}  \right)^2dt.
	\end{alignedat}
\end{equation}

	The key idea of the proof is to find some particular coefficients $\beta_n$ such that the ratio of \eqref{brasil} over \eqref{japao} is large, see \eqref{eq:phiM} and \eqref{eq:normalphiM}. 
	The choice of $\beta_n$ is such that only a finite number of coefficients $\beta_n$ is different of zero. Thus, let $M$ be a sufficient large integer and take 
\begin{equation*}
\beta_n=0\quad \forall n \not\in\{8M+k: \ 1\leq k\leq 8\}.
\end{equation*}
	Then, consider the initial conditions 
$$
	\varphi^{0,M}=\sum_M\beta_n\varphi_n
$$
	and its associated solution
 \begin{equation}\label{def_varphiM}
\varphi^M(\,\cdot\,,t)=\sum_M\alpha_n(t)\varphi_n\quad \forall  t \in (0,T),
\end{equation}
	where the symbol $\dis\sum_M$ stands for the sum extended to all indices $n$ of the form $n=8M+k$ with $1\leq k \leq 8$. The values of $\beta_n$, for 
	$n \in\{8M+k: \ 1\leq k\leq 8\}$, will be chosen in Section \ref{sec_above}.

\subsection{Estimate from below}\label{sec:3.2}

	In this section, let us use \eqref{brasil} to find an estimate like \eqref{eq:phiM}. To do this, let us begin with the inequality
\begin{equation*}\label{eq:frst}
    		\dis\sum_M\dfrac{1}{\lambda_n}\left(C_{1,n}e^{\mu_n^+T}+C_{2,n}e^{\mu_n^-T}\right)^2
    		\geq\dis\sum_M\dfrac{1}{\lambda_n}\left(\dfrac{3}{4} C_{2,n}^2e^{2\mu_n^-T}-3C_{1,n}^2e^{2\mu_n^+T}\right).
 \end{equation*}
 
 	First, let us assume that the constants $C_{1, 8M+k}$ and $\beta_{1, 8M+k}$ are bounded independently of $M$.  In fact, this will be proved in the next section, see Remark \ref{rmk:boundedC}.
	
  	Now, from \eqref{prop_lambda} and \eqref{def_mu}, we have that
\begin{equation}\label{eq:boundC1M}
C^2_{1, 8M+k}e^{2\mu^+_{ 8M+k}T}\leq Ce^{-CM^2T} \quad \forall k=1,\ldots, 8.
\end{equation}

On the other hand, using the notations
\begin{equation*}
(k-1/2)!=(k-1/2)(k-3/2)\cdots 1/2\quad \forall \ k\geq 1 \quad \hbox{and}\quad(-1/2)!=1,
\end{equation*}
we can expand the quotient $(a-\lambda_n+\sqrt{D_n})/\sqrt{D_n}$ in the definition of  $C_{2,n}$:
\begin{equation}\label{eq:C2}
\!\!\!\!\!\!\!\!
\begin{alignedat}{2}
\noalign{\smallskip}\dis
\dfrac{a-\lambda_n+\sqrt{D_n}}{\sqrt{D_n}}=\,&\left[\dfrac{2a}{\lambda_n+a}-\dfrac{2\lambda_n(\lambda_n-a)(a+b)}{(\lambda_n+a)^3}-\dfrac{\lambda_n-a}{\lambda_n+a}\sum_{k\geq2}\frac{(k-1/2)!}{k!}\left(\dfrac{4\lambda_n(a+b)}{(\lambda_n+a)^2}\right)^k\right]\\
\noalign{\smallskip}\dis
	=\,&\left[\dfrac{2a}{\lambda_n+a}-\dfrac{2\lambda_n(\lambda_n-a)(a+b)}{(\lambda_n+a)^3}- \dfrac{6\lambda_n^2(\lambda_n-a)(a+b)^2}{(\lambda_n+a)^5}+\mathcal{O}(\lambda_n^{-3})\right]
	\\
\noalign{\smallskip}\dis
	\approx~&\mathcal{O}(\lambda_n^{-1}),
\end{alignedat}
\end{equation}
	for $n$ large enough.

	This way, thanks to \eqref{prop_lambda}, we obtain
\begin{equation}\label{eq:boundC2M}
\inf_{ 1\leq k\leq 8}\left(\dfrac{a-\lambda_{2, 8M+k}+\sqrt{D_{2, 8M+k}}}{\sqrt{D_{2, 8M+k}}}\right)^2\geq \dfrac{C}{M^4},
\end{equation}
for $M$ large enough and for some positive constant $C$ independent of $M$. 

	Finally, combining  \eqref{brasil} with \eqref{prop_lambda},  \eqref{eq:boundC1M}, \eqref{eq:boundC2M} and the convergence $\mu^-_n\to -(a+b)$, one has:
\begin{equation}\label{roma}
\|\varphi^M(\,\cdot\,,0)\|^2\geq \dfrac{C_0}{M^6},
\end{equation}
	 for $M$ large enough and for some positive constant $C_0$ independent of $M$. 

\subsection{Estimate from above}\label{sec_above}

	In this section, we are going to obtain an estimate for \eqref{japao}. Let us first set
\begin{equation}\label{A_1}
A_1:= 96\pi R^2\int_0^T\left(\sum_M\gamma_n C_{1,n}e^{(a+b+\mu_n^+)(T-t)}   \right)^2dt
\end{equation}
and 
\begin{equation}\label{A_2}
A_2:=96\pi R^2\int_0^T\left(\sum_M\gamma_nC_{2,n}e^{(a+b+\mu_n^-)(T-t)}  \right)^2dt.
\end{equation}
We will analyze the estimates for \eqref{A_1} and \eqref{A_2}, separately.

\

\begin{claim} There exists a constant $C>0$, independent of $M$, such that 
\begin{equation}\label{eq:A_1}
	A_1\leq {C\over M^{10}},
\end{equation}
	for $M$ large enough.
\end{claim}

\

\begin{proof}
	Let us begin using \eqref{def_mu} to split the term:
\begin{equation*}
e^{(a+b+\mu_n^+)(T-t)}=e^{(a+2b-\lambda_n)(T-t)}e^{B_n(T-t))}
\end{equation*}
where $B_n:=-\mu_n^--a-b\to 0$ as $n\to \infty$. Also, from \eqref{prop_lambda}, we have
\begin{equation*}
	\begin{array}{lll}
		e^{(a+2b-\lambda_{8M+k})(T-t)}=e^{\left[a+2b-\frac{\pi^2}{R^2}\left(8M+\frac{1}{2}\right)^2\right](T-t)}
		e^{\left[-\frac{\pi^2}{R^2}(16Mk+k+k^2)+\varepsilon_{8M+k}\right](T-t)}.
	\end{array}
\end{equation*}

	Now, let us rewrite $A_1$ as follows:
\begin{equation*}
	A_1=96\pi R^2\int_0^Te^{(2a+4b-2\frac{\pi^2}{R^2}(8M+\frac{1}{2})^2)(T-t)}g_M(t) \, dt,
\end{equation*}
where $g_M(t):=f_M(t)^2$ with $f_M$ given by
\begin{equation*}
f_M(t):=\sum_{k=1}^8\gamma_{8M+k}C_{1,8M+k}e^{ \left[   -\frac{\pi^2}{R^2}(16Mk+k+k^2)+\varepsilon_{8M+k} +B_{8M+k}    \right](T-t)}.
\end{equation*}

	Let us now integrate by parts ten times the integral term in $A_1$. Then, we have:
\begin{equation}\label{eq:eq1}
	\begin{array}{lcl}
		\dis\int_0^Te^{(2a+4b-\frac{2\pi^2}{R^2}(8M+\frac{1}{2})^2)(T-t)}g_M(t) \ dt
		\!&=&\!\dis\sum_{j=0}^9\dfrac{e^{(2a+4b-\frac{2\pi^2}{R^2}(8M+\frac{1}{2})^2)T} g_M^{(j)}(0)-g_M^{(j)}(T)}{(2a+4b-\frac{2\pi^2}{R^2}(8M+\frac{1}{2})^2)^{j+1}}\\
		\noalign{\smallskip}\dis
		&&\!\!\dis+\int_0^T\dfrac{e^{(2a+4b-\frac{2\pi^2}{R^2}(8M+\frac{1}{2})^2)(T-t)}}{(2a+4b-\frac{2\pi^2}{R^2}(8M+\frac{1}{2})^2)^{10}}g_M^{(10)}(t) \ dt.
	\end{array}
\end{equation}

	Since  the constants $C_{1, 8M+k}$ will be chosen to be bounded independently of $M$ and  $\varepsilon_{8M+k}$,  $B_{8M+k}$ and  $\gamma_{8M+k}$ 
	are bounded independently of $M$, 
	we have that $|f_M^{(j)}|=\mathcal{O}(M^j)$ and $g^{(j)}_M=\mathcal{O}(M^j)$ for all $j\geq1$ and  $M$ large enough. 
	Then, 
\begin{equation*}
	\begin{array}{lcl}
		\!\dis\sum_{j=0}^9\dfrac{g_M^{(j)}(T)}{(2a+4b-\frac{2\pi^2}{R^2}(8M+\frac{1}{2})^2)^{j+1}}=\mathcal{O}(M^{-2}).
	\end{array}
\end{equation*}
	Therefore, in order to obtain \eqref{eq:A_1} we need to impose conditions on $g_M^{(j)}(T)$:
\begin{equation}\label{cond_g}
g_M^{(0)}(T)=g_M^{(1)}(T)=\cdots=g_M^{(8)}(T)=g_M^{(9)}(T)=0.
\end{equation}

	Notice that these conditions are fulfilled if the constants $C_{1,8M+k}$ $(1\leq k\leq 8)$ satisfy five linear equations corresponding to the identities $f^{(0)}_M(T) = f_M^{(1)} (T) =f_M^{(2)}(T) = f_M^{(3)}(T)=f_M^{(4)}(T) =   0$. More precisely, the constants $C_{1,8M+k}$ $(1\leq k\leq 8)$ should satisfy:
\begin{equation}\label{sistema_1}
\left\{
    \begin{array}{l}\dis
    		\sum_{k=1}^8\gamma_{8M+k}  \left(  -\frac{\pi^2}{R^2}(16Mk+k+k^2)+\varepsilon_{8M+k} +B_{8M+k}  \right)^jC_{1,8M+k}=0,\\
    		  \noalign{\smallskip}\dis
    		\hbox{for } j=0,1,2,3,4.
    \end{array}
    \right.
\end{equation} 
\begin{remark}
	In the linear system \eqref{sistema_1}, we have five linear equations and eight unknowns $C_{1,8M+k}$. Hence, since \eqref{sistema_1} is a linear homogeneous system,
	 the space of solution has, at least, dimension 1. Therefore, it is not difficult to choose a nontrivial solution to \eqref{sistema_1} bounded independently of $M$.
\end{remark}

	Finally, using \eqref{eq:eq1}, \eqref{cond_g} and the following bounds
\begin{equation*}
\dfrac{e^{(2a+4b-\frac{2\pi^2}{R^2}(8M+\frac{1}{2})^2)T}}{(2a+4b-\frac{2\pi^2}{R^2}(8M+\frac{1}{2})^2)^{j+1}}|g^{(j)}_M(0)|\leq C e^{-CM^2}\dfrac{1}{M^{j+2}}<\dfrac{C}{M^{10}}
\end{equation*}
for  $0\leq j\leq 9$ and
\begin{equation*}
\left|\int_0^T\dfrac{e^{(2a+4b-\frac{2\pi^2}{R^2}(8M+\frac{1}{2})^2)(T-t)}}{(2a+4b-\frac{2\pi^2}{R^2}(8M+\frac{1}{2})^2)^{10}}g_M^{(10)}(t) \ dt\right|\leq \int_0^T\dfrac{1}{(CM)^{20}}CM^{10} \ dt=\dfrac{C}{M^{10}},
\end{equation*}
	for $M$ large enough, we deduce \eqref{eq:A_1}.

\end{proof}

\

\begin{claim} There exists a constant $C>0$, independent of $M$, such that 
\begin{equation}\label{eq:A2}
	A_2\leq {C\over M^{12}},
\end{equation}
	for $M$ large enough.
\end{claim}

\

\begin{proof} First, let us rewrite and expand $\mu^-_n$ in the following way:
\begin{equation*}
\begin{alignedat}{2}
		\mu^-_n=~&\dfrac{\lambda_n+a}{2}\left(-1+\sqrt{1-\dfrac{4\lambda_m(a+b)}{(\lambda_n+a)^2}}\right)	\\
		\noalign{\smallskip}\dis
		=~&\dis-\dfrac{\lambda_n+a}{4}\sum_{k\geq1}\dfrac{(k-3/2)!}{k!}\left[\dfrac{4\lambda_n(a+b)}{(\lambda_n+a)^2}\right]^k.
\end{alignedat}
\end{equation*}
	
	Then, the exponent in the expression of $A_2$ can be split in the form: 
\begin{equation*}
		e^{(a+b+\mu^-_n)(T-t)}=e^{{a(a+b)\over\lambda_n+a}(T-t)}e^{Y_n(T-t)},
\end{equation*}
	where
\begin{equation*}
Y_n:=-\dfrac{\lambda_n+a}{4}\sum_{k\geq2}\dfrac{(k-3/2)!}{k!}\left[\dfrac{4\lambda_n(a+b)}{(\lambda_n+a)^2}\right]^k.
\end{equation*}	
	
	Since $e^x=1+x+\mathcal{O}(x^2)$ for $|x|<1$,  we see that
\begin{equation}\label{peru}
e^{{a(a+b)\over\lambda_n+a}(T-t)}=1+ \dfrac{a(a+b)}{\lambda_n+a}(T-t)+\mathcal{O}(\lambda_n^{-2}),
\end{equation}
	for $n$ large enough.

	Now, since $\mu_n^-\to -(a+b)$, we have that
\begin{equation*}
|Y_n(T-t)|=\left|\left(a+b+\mu^-_n-\dfrac{a(a+b)}{\lambda_n+a}\right)(T-t)\right|<1,
\end{equation*}
and
\begin{equation}\label{argentina}
e^{Y_n(T-t)}=1-\dfrac{\lambda_n^2(a+b)^2}{(\lambda_n+a)^3}(T-t)+\mathcal{O}(\lambda_n^{-2}),
\end{equation}
	where we have used that $Y_n=-{\lambda_n^2(a+b)^2\over(\lambda_n+a)^3}+\mathcal{O}(\lambda_n^{-2})$,  for $n$ large enough.

	Therefore, from \eqref{peru} and \eqref{argentina}, we have that
\begin{equation}\label{eq:expmu}
e^{(a+b+\mu^-_n)(T-t)}=1-\dfrac{\lambda_n^2(a+b)^2}{(\lambda_n+a)^3}(T-t)+ \dfrac{a(a+b)}{\lambda_n+a}(T-t)+\mathcal{O}(\lambda_n^{-2}).
\end{equation}


	On the other hand, using \eqref{eq:C2} and \eqref{eq:expmu}, we obtain
\begin{equation}\label{belgica}
\begin{alignedat}{2}
	\gamma_nC_{2,n}e^{(a+b+\mu^-_n)(T-t)}=~&\gamma_n\dfrac{\beta_n}{2}\bigg[\left(\dfrac{2a}{\lambda_n+a}   -\dfrac{2\lambda_n(\lambda_n-a)(a+b)}{(\lambda_n+a)^3}- \dfrac{6\lambda_n^2(\lambda_n-a)(a+b)^2}{(\lambda_n+a)^5}\right)\\
	\noalign{\smallskip}\dis
	&+(T-t)\bigg(-\dfrac{2\lambda_n^2 (a+b)^2a}{(\lambda_n+a)^4}   + \dfrac{2\lambda_n^3(\lambda_n-a)( a+b)^3}{(\lambda_n+a)^6} +     \dfrac{2a^2( a+b)}{(\lambda_n+a)^2}   \\
	\noalign{\smallskip}\dis
	&-\dfrac{2\lambda_n(\lambda_n-a)a( a+b)^2}{(\lambda_n+a)^4   }\bigg)+     \mathcal{O}(\lambda_n^{-3})\bigg],
\end{alignedat}
\end{equation}
	for $n$ large enough.
	
	In order to obtain \eqref{eq:A2}, from the previous identity, we should impose that 
\begin{equation}\label{sistema_2}
\sum_M\gamma_n\left(\dfrac{a}{\lambda_n+a}-\dfrac{\lambda_n(\lambda_n-a)(a+b)}{(\lambda_n+a)^3}- \dfrac{3\lambda_n^2(\lambda_n-a)(a+b)^2}{(\lambda_n+a)^5}\right)\beta_n=0
\end{equation}
and 
\begin{equation}\label{sistema_3}
\sum_M\gamma_n\left(\dfrac{\lambda_n^2 (a+b)^2a}{(\lambda_n+a)^4}   - \dfrac{\lambda_n^3(\lambda_n-a)( a+b)^3}{(\lambda_n+a)^6} -     \dfrac{a^2( a+b)}{(\lambda_n+a)^2} +\dfrac{\lambda_n(\lambda_n-a)a( a+b)^2}{(\lambda_n+a)^4} \right)\beta_n=0.
\end{equation}
\begin{remark}\label{rmk:boundedC}
	The expression \eqref{def_C_n}, links the choices of $C_{1,n}$ and $\beta_n$. This way, \eqref{sistema_1}, \eqref{sistema_2} and \eqref{sistema_3} is a linear homogeneous system of seven equations 
	and eight unknowns. As before, the constants $C_{1, 8M+k}$ (and equivalently $\beta_{8M+k}$) can be chosen to be bounded independently of $M$.
\end{remark}

	 Finally, from \eqref{belgica}, \eqref{sistema_2} and  \eqref{sistema_3}, we have that there exists a constant $C$, independent of $M$, such that
\begin{equation*}
	\begin{array}{lll}
		\dis \left|\sum_M \gamma_nC_{2,n}e^{(a+b+\mu_n^-)(T-t)}\right|
		&\leq&\dfrac{C}{ M^6},
	\end{array}
\end{equation*}
	for $M$ large enough, and this leads to \eqref{eq:A2}.
\end{proof}

\

	As a consequence of the estimates of \eqref{eq:A_1} and \eqref{eq:A2}, we deduce from \eqref{japao} that
\begin{equation}\label{estimate_normal}
\iint_\Sigma e^{2(a+b)(T-t)}\left|\dfrac{\partial\varphi^M}{\partial n}\right|^2d\Gamma dt\leq  \dfrac{C}{M^{10}},
\end{equation}
	for $M$ large enough.

\subsection{Construction of non-controllable initial data}\label{contruction}

	From the results obtained in Sections \ref{sec:3.1}, \ref{sec:3.2} and \ref{sec_above}, it is clear that there is no constant $C$ such that \eqref{des-obs} holds.  
	Consequently, \eqref{eq_stokes+memory} is not null-controllable.

	For the sake  of completeness, let us construct explicitly initial data $\overline{y}_0 \in H(\Omega)$ such that, for all $v\in L^2(\gamma\times(0,T))$, the corresponding state does 
	not vanish at $t=T$.

	First, note from \eqref{brasil} and \eqref{roma} that, for each $M$ large enough, there exists $k_0$ such that  $1\leq k_0\leq 8$ and
\begin{equation}\label{herika}
\|\varphi_{ 8M+k_0}\|^2\left(  C_{1, 8M+k_0}e^{\mu ^+_{ 8M+8k_0}T}  + C_{2, 8M+k_0}e^{\mu^-_{8M+k_0}T} \right)^2\geq \dfrac{C_0}{ 8M^{6}}.
\end{equation}

	Then, let us define
\begin{equation}\label{def_y^0}
\overline{y}_0=\sum_{\ell\geq 1}\dfrac{1}{\ell^{3/4}}\dfrac{\varphi_{ 8\ell+k_0}}{\|\varphi_{ 8\ell+k_0}\|}.
\end{equation}
	It is not difficult to see that  $\overline{y}_0\in H(\Omega)$. 
	
	Let us now prove that $\overline{y}_0$ cannot be steered to the rest. Indeed, arguing by contradiction, let $v\in L^2(\Sigma)$ be such that the solution to (\ref{eq_stokes+memory}) satisfies  $y(\cdot\,,T)=0$. 
	Then, in particular, we have that 
	\begin{equation}\label{dualidade}
		\!\!\!\int_\Omega \overline{y}_0(x)\varphi^{M}(x,0) \ dx\dis=\!\!\iint_\Sigma v 			\dfrac{\partial \varphi^{M}}{\partial n}\ d\Gamma \ dt
		 + \ b\!\int_0^T \!\!\! 			\int_0^te^{-a(t-s)}\left( \int_{\partial\Omega} v(\sigma,s)\dfrac{\partial 						\varphi^{M}}{\partial n}(\sigma,t) \ d\Gamma \right)ds \ dt,
	\end{equation}
	where $\varphi^{M}$ is defined in \eqref{def_varphiM}.

	Using \eqref{def_y^0} and the orthogonality of $\{\varphi_n\}_{n\in \mathbb{N}}$, we obtain
\begin{equation*}
	\begin{alignedat}{2}
		\dis\int_\Omega \overline{y}_0(x)\varphi^{M}(x,0) \ dx
		=&\dfrac{1}{M^{3/4}}\|\varphi_{ 8M+k_0}\|\left( C_{1, 8M+k_0}e^{\mu ^+_{ 8M+k_0}T}  + C_{2, 8M+k_0}e^{\mu^-_{ 8M+k_0}T}   \right)
		\end{alignedat}
\end{equation*}
and, in view of \eqref{herika}, we find
\begin{equation}\label{mexico}
\left|\int_\Omega \overline{y}_0(x)\varphi^{M}(x,0) \ dx\right|\geq\dfrac{C_1}{M^{15/4}},
\end{equation}
	for $M$ large enough and for some positive constant $C_1$ independent of $M$. 

	 On the other hand, taking into account \eqref{estimate_normal}, we see that the other terms in \eqref{dualidade} can be bounded as follows
 \begin{equation}\label{marrocos}
 	\left| \iint_\Sigma v(\sigma,t) \dfrac{\partial \varphi^{M}}{\partial n}(\sigma,t) \ d\Gamma  \ dt  \right|
	\leq \|v\|_{L^2(\Sigma)}\left\|   \dfrac{\partial \varphi^{M}}{\partial n}\right\|_{L^2(\Sigma)}\leq \dfrac{C_2}{M^5}
 \end{equation}
 	and
 \begin{equation}\label{australia}
	\left|\int_0^T\int_0^te^{-a(t-s)}\left( \int_{\partial\Omega} v(\sigma,s)\dfrac{\partial \varphi^{M}}{\partial n}(\sigma,t) d\Gamma \right) \ ds \ dt\right| 
	\leq C\|v\|_{L^2(\Sigma)}\left\|   \dfrac{\partial \varphi^{M}}{\partial n}\right\|_{L^2(\Sigma)}\leq \dfrac{C_3}{M^5},
 \end{equation} 
 	for $M$ large enough and for some positive constants $C_2$ and $C_3$ independent of $M$. 
 
 	Consequently, \eqref{mexico}, \eqref{marrocos} and \eqref{australia} lead to
 \begin{equation*}
 \dfrac{C_1}{M^{15/4}}\leq \dfrac{C_4}{M^5},
 \end{equation*}	
	 which is an absurdity.

%
%


\section{Additional comments and questions}\label{coments}

\subsection{Lack of null controllability for the two-dimensional Stokes equations with a memory term}

	In this section, let us present the key points, similar to those in Section \ref{section_eigenfunction}, to obtain lack of null controllability for \eqref{eq_stokes+memory} in the two-dimensional case.
	For simplicity, the coordinates of a generic point in $\Omega$ will be denoted by $x$ and $y$, where $\Omega$ is a ball of radius $R$ centered at the origin.

	As in Section  \ref{section_eigenfunction}, we compute $(\varphi,q$) and $\lambda$ satisfying  \eqref{prob_autovalor}.  Then, considering $\varphi=\nabla\times \psi$ and setting $w=\nabla\times\varphi$,
	we get to the following eigenvalue problem
\begin{equation}\label{italia}
\left\{
    \begin{array}{lll}
    		\dis -r^2w^{\prime\prime}-rw^{\prime}= \lambda r^2w& \hbox{in} & (0,R),\\
    		\noalign{\smallskip}\dis
    		-r^2\psi^{\prime\prime}-r\psi^{\prime}= r^2w& \hbox{in} & (0,R),\\
    		\noalign{\smallskip}\dis
    		\psi(R)=0, \ \ \psi^{\prime}(R)=0, \  \ \lambda>0.&  & 
    \end{array}
\right.
\end{equation}
	Here, we are using the notation $r=\sqrt{x^2+y^2}$ for any $(x,y)\in \Omega$.
	
	The solutions to \eqref{italia}$_1$ are linear combinations of Bessel functions, i.e.  $w(r)=c_1J_0(\lambda^{1/2}r)+c_2Y_0(\lambda^{1/2}r)$, where $J_0$ and $Y_0$ 
	are Bessel functions of zero order and of first and second kind, respectively. Since $Y_0(0)=-\infty$, we can take $c_2=0$ and $w(r)=J_0(\lambda^{1/2}r)$.
	Then, from \eqref{italia}$_2$, we have that $$\left(r\psi^{\prime}\right)^{\prime}=-rw.$$ Now, integrating twice with respect to $r$ and using the boundary conditions  \eqref{italia}$_3$, 
	we obtain
\begin{equation*}
	\psi(r)=-\int_r^R\dfrac{1}{\sigma}\left(\int_\sigma^RsJ_0(\lambda^{1/2}s)\,ds\right)\,d\sigma.
\end{equation*}

	In view of the identity $\dfrac{d}{dr}\left[ r^pJ_p(r)   \right]=r^pJ_{p-1}(r)$, we have:
\begin{equation*}
	\int_\sigma^RsJ_0(\lambda^{1/2}s)\,ds=\dfrac{1}{\lambda^{1/2}}\left[ RJ_1(\lambda^{1/2}R)  -\sigma J_1(\lambda^{1/2}\sigma)    \right].
\end{equation*}
	Therefore, we deduce
\[
			\psi(r)
			=-\dfrac{RJ_1(\lambda^{1/2}R)}{\lambda^{1/2}}(\log R-\log r)+\dfrac{1}{\lambda}\int_{\lambda^{1/2}r}^{\lambda^{1/2}R}J_1(\sigma)\,d\sigma.
\]
	
	If we choose $\lambda$ such that $J_1(\lambda^{1/2}R)=0$, we obtain that
\begin{equation*}
	\psi(r)=\dfrac{1}{\lambda}\int_{\lambda^{1/2}r}^{\lambda^{1/2}R}J_1(\sigma)\,d\sigma.
\end{equation*}
	Then, we can define the following eigenvalues with the corresponding eigenfunctions:
\begin{equation}\label{eq:bessel}
\left\{
\begin{array}{l}
\lambda_n^{1/2}R=j_{1,n}\\
\noalign{\smallskip}\dis
\psi_n(r)=\dfrac{1}{\lambda_n}\int_{\lambda_n^{1/2}r}^{\lambda^{1/2}_nR}J_1(\sigma)\,d\sigma\\
\noalign{\smallskip}\dis
q_n\equiv0
\\
\noalign{\smallskip}\dis
\varphi_n(x,y)=\dfrac{J_1(\lambda^{1/2}_nr)}{\lambda^{1/2}_nr}\left(   -y,x\right),
\end{array}
\right.
\end{equation}
where $j_{1,n}$ is the {\it n}-th positive root of $J_1$.

	Thanks to \cite[Lemma 1]{Lorch}, $\lambda_n$ satisfies the following inequality:
\begin{equation}\label{est_lambda_case_2D}
\dfrac{\pi^2}{R^2}\left(n+\dfrac{1}{8}\right)^2\leq \lambda_n\leq \dfrac{\pi^2}{R^2}\left(n+\dfrac{1}{4}\right)^2    \quad    \forall\, n\geq 1.
\end{equation}
{
	Taking into account \eqref{eq:bessel}$_1$, a simple computation gives us:
%
\begin{equation}\label{eq:normalbessel}
	 \dfrac{\partial \varphi_n}{\partial n}\biggr|_{\partial\Omega}=J_1^{\prime}(\lambda_n^{1/2}R)\left(-\dfrac{y}{R},\dfrac{x}{R}\right).
\end{equation}

	On the other hand, thanks to the inequality \eqref{est_lambda_case_2D}, we also obtain the following bound:
\begin{equation}\label{est_varphi_2D_case}
\begin{array}{lll}
    		\dis \|\varphi_n\|^2
    		&=&\dis\dfrac{1}{\lambda_n}\dis \int_\Omega [J_1(\lambda^{1/2}_nr)]^2\, dx\,dy \\
    		\noalign{\smallskip}\dis
    		&=& \dfrac{2\pi}{\lambda^2_n}\dis \int_0^{j_{1,n}} [J_1(s)]^2s \,ds\\
    		\noalign{\smallskip}\dis
    		&\geq& \dfrac{2\pi}{\lambda^2_n}\dis \int_0^{1} J_1^2(r)r \,dr\\
    		\noalign{\smallskip}\dis
		&\geq&\dfrac{2\pi C}{\lambda^2_n}.
    \end{array}
\end{equation}

	Finally, analogously to the three-dimensional case, we can define $\gamma_n:=J_1^{\prime}(\lambda_n^{1/2}R)$. Thanks to \eqref{eq:bessel}$_1$, it is not difficult to see that
	$\gamma_n=J_0(\lambda_n^{1/2}R)$ and therefore it is bounded independent of $n$.

	In view of \eqref{est_lambda_case_2D}, \eqref{eq:normalbessel}, \eqref{est_varphi_2D_case} and the boundedness of $\gamma_n$,
	we can adapt the proof of Theorem \ref{theo_main} to obtain same result in the two-dimensional case. 
	
\begin{remark}
	The main difference between the two and three dimensional cases is the fact that estimate \eqref{roma} will be slightly different, i.e. in the two-dimensional case we obtain  
	a bound of order $ \mathcal{O}(M^{-8})$. 
	On the other hand, the estimate \eqref{estimate_normal} holds for both cases and it is sufficient to conclude the proof of Theorem  \ref{theo_main} in the two-dimensional case. 
\end{remark}

}

\subsection {Heat equation with memory}

	The non-null controllability results obtained in \cite{Guerrero}, for equation  \eqref{eq_guerrero}, can be extended to more general situations. More precisely, using 
	similar arguments and computations done in the previous sections, it is not difficult to obtain a similar result for a heat equation with heat flux memory:
\begin{equation*}\label{eq_heatDE}
\left\{
    \begin{array}{lll}\dis
    		y_t-\Delta y-b\int_0^te^{-a(t-s)} \Delta y(\cdot\,,s) \, ds = 0 & \hbox{in} & Q,\\
    		y=v& \hbox{on} & \Sigma,\\
    		\noalign{\smallskip}\dis
    		y(\cdot\,,0)=y^0 &\hbox{in} &\Omega.
    \end{array}
    \right.
\end{equation*}


	It would be interesting to investigate which are the optimal conditions on a time-dependent memory kernel $K$ such that Theorem \ref{theo_main} still holds for the systems like
\begin{equation*}
\left\{
    \begin{array}{lcl}\dis
    		y_t-\Delta y-\int_0^tK(t-s) \Delta y(\cdot\,,s) \ ds = 0 & \hbox{in} & Q,\\
    		y=v& \hbox{on} & \Sigma,\\
    		\noalign{\smallskip}\dis
    		y(\cdot\,,0)=y^0 &\hbox{in} &\Omega.
    \end{array}
    \right.
\end{equation*}
	Some results in the one-dimensional case were obtained in \cite{Halanay}.

\subsection{Controls with less components}

	The approximate controllability result for \eqref{eq_stokes+memory_distributed}, proved in \cite{Doubova}, can be improved with respect to the number of scalar controls. 
	More precisely, we will use similar ideas from \cite{Lions_Zuazua} to show a new unique continuation property for the solutions to \eqref{1a}, i.e.: 
\begin{equation*}
\varphi_i=0 \quad \hbox{in}\quad \omega\times(0,T), \quad\hbox{for} \quad 1\leq i<N\quad \Rightarrow\quad\varphi=0\quad\Omega\times(0,T),
\end{equation*}
	where $\omega\subset \Omega$ is an open subset.
	Indeed, let us first notice that the free divergence condition for $\varphi$ implies that $p$ is harmonic in $\Omega$ with respect to $x$ for all $t\in(0,T)$.
	On the other  hand, since $\varphi_i=0$ in $\omega\times(0,T)$, for $1\leq i<N$, we deduce that 
\begin{equation*}
\dfrac{\partial p}{\partial x_i}=0\quad\hbox{in}\quad\omega\times(0,T) \quad\hbox{for} \quad 1\leq i<N.
\end{equation*}
	
	Hence, elliptic unique continuation guarantees that
\begin{equation*}
\dfrac{\partial p}{\partial x_i}=0\quad\hbox{in}\quad\Omega\times(0,T) \quad\hbox{for} \quad1\leq  i<N
\end{equation*}
	and then $p$ is a function that only depends of $x_N$.

	Now, one can see that $\varphi_i$, for $1\leq i<N$, satisfies the following backward heat equation with memory:
\begin{equation*}
\left\{
    \begin{array}{lll}\dis
    		-\varphi_{i,t}-\Delta \varphi_i-b\int_t^Te^{-a(s-t)}\Delta \varphi_i(\cdot\,,s) \ ds = 0 & \hbox{in} & Q,\\
    		\varphi_i=0& \hbox{on} & \Sigma,\\
    		\noalign{\smallskip}\dis
    		\varphi_i(\cdot\,, T)=\varphi_i^0 &\hbox{in} &\Omega.
    \end{array}
    \right.
\end{equation*} 
Since the unique continuation result \cite[Lemma 2.3]{Doubova} holds for both heat equation and Stokes equations with memory and the fact that $\varphi_i=0$ in $\omega\times(0,T)$, for $1\leq i<N$, 
we have
\begin{equation*}
\varphi_i=0\quad\hbox{in}\quad\Omega\times(0,T)\quad\hbox{for}\quad 1\leq i<N.
\end{equation*}

	Finally, the free divergence condition and the Dirichlet boundary condition lead to
\begin{equation*}
\varphi=0\quad\hbox{in}\quad \Omega\times(0,T).
\end{equation*}

\begin{remark}
	In general, in the three-dimensional case, the approximate controllability result does not hold using only one control. Indeed, let us consider $L>0$ and $\Omega=G\times (0,L)$, where
	$G\subset \mathbb{R}^2$ is a bounded domain such that there exists an eigenfunction of the Stokes equations in $G$ in which the corresponding pressure is zero (for example, if $G$ is 
	a ball then this property holds). More precisely, there exist a nontrivial vector field $u=(u_1,u_2)$ and a number $\lambda>0$ such that 
\begin{equation*}
\left\{
    \begin{array}{lll}\dis
    		-\Delta_{(x,y)} u=\lambda u & \hbox{in} & G,\\
    		\noalign{\smallskip}\dis
    		\nabla_{(x,y)} \cdot u=0 &\hbox{in} &G,\\
    		\noalign{\smallskip}\dis
    		u=0& \hbox{on} & \partial G.
    \end{array}
    \right.
\end{equation*}

	To sum up, let us define 
$\varphi(x,y,z,t):=\alpha(t)\sin\left({\pi z}/{L}\right)\left(u_1(x,y),u_2(x,y),0\right)$, where $\alpha$ is the solution to \eqref{EDO_alpha} with $\beta_n=1$ and $\lambda_n=\lambda+\pi^2/L^2$.
Then,  $\varphi$ is a nonzero solution to \eqref{1a} such that $\varphi_3\equiv0$.

\end{remark}

\subsection{Hyperbolic equations with memory}

	Differently to the case of the heat and Stokes equations with memory, the exact controllability for wave equation with memory holds. More precisely, for a hyperbolic integro-differential equation 
	of the form
\begin{equation*}
\left\{
    \begin{array}{lll}\dis
    		y_{tt}-a(t)\Delta y+b(t)y_t+c(t)y-\int_0^t K(t,s)\Delta y(\cdot\,,s) \ ds=0 & \hbox{in} & \Omega\times (0,T),\\
    		y=v& \hbox{on} & \partial \Omega\times (0,T),\\
    		\noalign{\smallskip}\dis
    		y( \cdot\,,0)=0, \quad y_t(\cdot, 0)=0 &\hbox{in} &\Omega,
    \end{array}
    \right.
\end{equation*}
	the exact controllability holds as well as the kernel $K=K(t,s)$ belongs to $C^2(\mathbb{R}^2_+)$. For more details, see \cite{Kim}.

	It would be interesting analyze if the exact controllability results obtained in \cite{Kim} can be extended to the hyperbolic Stokes equation with memory:
\begin{equation*}\label{hypeq_stokes+memory}
\left\{
    \begin{array}{lcl}\dis
    		y_{tt}-\Delta y-\int_0^tK(t,s)\Delta y(\cdot\,,s) \ ds +\nabla p = 0 & \hbox{in} & Q,\\
    		  \noalign{\smallskip}\dis
    		\nabla\cdot y=0& \hbox{in} &Q,\\
    		y=v& \hbox{on} & \Sigma,\\
    		\noalign{\smallskip}\dis
    		y(\cdot\,,0)=0,\quad y_t(\cdot\,,0)=0 &\hbox{in} &\Omega.
    \end{array}
    \right.
\end{equation*}

\subsection{Nonlinear systems with memory}
	
	Recall that the null and approximate controllability of \eqref{old_2} and \eqref{old_1} are open questions. It would be very interesting to see whether or not the effect of the nonlinear terms 
	is sufficient to modify the controllability properties of the linearized systems. This is the case, for instance, for the equation studied in \cite{coron-lissy}.

\

\noindent\textbf{Acknowledgments:} This work has been partially done while the second author was visiting the 
		Universidad de Sevilla (Seville, Spain). He wishes to thank the members of the IMUS (Instituto de Matem\'aticas de la Universidad de Sevilla) for their kind hospitality.

\end{document}